\documentclass[a4paper,11pt]{amsart}
\usepackage{graphicx}
\usepackage{mathptmx}
\usepackage{amsmath}
\usepackage{amssymb}
\usepackage{enumitem}
\usepackage{xcolor}
\usepackage{xparse}
\NewDocumentCommand{\eulerian}{omm}
 {%
  \genfrac<>{0pt}{}{#2}{#3}%
  \IfValueT{#1}{_{\!#1}}%
 }

 \newmuskip\pFqmuskip

\newcommand*\pFq[6][8]{%
  \begingroup 
  \pFqmuskip=#1mu\relax
  \mathchardef\normalcomma=\mathcode`,
  \mathcode`\,=\string"8000
  \begingroup\lccode`\~=`\,
  \lowercase{\endgroup\let~}\pFqcomma
  {}_{#2}F_{#3}{\left(\genfrac..{0pt}{}{#4}{#5}\bigg|#6\right)}%
  \endgroup
}
\newcommand{\pFqcomma}{{\normalcomma}\mskip\pFqmuskip}

\newtheorem{theorem}{Theorem}

\newtheorem{corollary}[theorem]{Corollary}

\begin{document}

\title[A Study on $q$-analogues of Catalan-Daehee numbers and polynomials]{A Study on $q$-analogues of Catalan-Daehee numbers and polynomials}

\author{Yuankui Ma}
\address{School of Science, Xi’an Technological University, Xi’an, 710021, Shaanxi, P. R. China}
\email{mayuankui@xatu.edu.cn}

\author{Taekyun  Kim}
\address{School of Science, Xi’an Technological University, Xi’an, 710021, Shaanxi, P. R. China\\
Department of Mathematics, Kwangwoon University, Seoul 139-701, Republic of Korea}
\email{tkkim@kw.ac.kr}

\author{DAE SAN KIM}
\address{Department of Mathematics, Sogang University, Seoul 121-742, Republic of Korea}
\email{dskim@sogang.ac.kr}

\author{Hyunseok Lee}
\address{Department of Mathematics, Kwangwoon University, Seoul 139-701, Republic of Korea}
\email{luciasconstant@kw.ac.kr}


\maketitle

\begin{abstract}

Catalan-Daehee numbers and polynomials, generating functions of which can be expressed as  $p$-adic Volkenborn integrals on $\mathbb{Z}_p$, were studied previously. The aim of this paper is to introduce $q$-analogues of the Catalan-Daehee numbers and polynomials with the help of $p$-adic $q$-integrals on $\mathbb{Z}_{p}$. We derive, among other things, some explicit expressions for the $q$-analogues of the Catalan-Daehee numbers and polynomials. 
\end{abstract}

\section{Introduction and preliminaries}
In recent years, many special numbers and polynomials have been studied by using several different tools such as combinatorial methods, generating functions, $p$-adic analysis, umbral calculus, differential equations, probability theory, special functions and analytic number theory. Catalan-Daehee numbers and polynomials were studied in [10] and several properties and identities associated with those numbers and polynomials were derived by utilizing umbral calculus techniques. The family of linear differential equations arising from the generating function of Catalan–-Daehee numbers were considered in [11] in order to derive some explicit identities involving Catalan-–Daehee numbers and Catalan numbers. In [6], $w$-Catalan polynomials were introduced as a generalization of Catalan polynomials and many symmetric identities in three variables related to the $w$-Catalan polynomials and analogues of alternating power sums were obtained by means of $p$-adic fermionic integrals.
The aim of this paper is to introduce $q$-analogues of the Catalan-Daehee numbers and polynomials with the help of $p$-adic $q$-integrals on $\mathbb{Z}_p$, and derive some explicit expressions and identities related to those numbers and polynomials. For the rest of this section, we recall the necessary facts that are needed throughout this paper. \par
Let $p$ be a fixed odd prime number. Throughout this paper, $\mathbb{Z}_{p},\mathbb{Q}_{p}$ and $\mathbb{C}_{p}$ denote respectively the ring of $p$-adic integers, the field of $p$-adic rational numbers and the completion of the algebraic closure of $\mathbb{Q}_{p}$. The $p$-adic norm $|\cdot|_{p}$ is normalized as $|p|_{p}=\frac{1}{p}$. Let $q$ be an indeterminate in $\mathbb{C}_{p}$ with $|1-q|_{p}<p^{-\frac{1}{p-1}}$. The $q$-analogue of $x$ is defined by $\displaystyle [x]_{q}=\frac{1-q^{x}}{1-q}\displaystyle$. Note that $\displaystyle \lim_{q\rightarrow 1}[x]_{q}=x.\displaystyle$ \par 
Let $f$ be a uniformly differentiable function on $\mathbb{Z}_{p}$. Then the $p$-adic $q$-integral on $\mathbb{Z}_{p}$ is defined by Kim as 
\begin{align}
\int_{\mathbb{Z}_{p}}f(x)d\mu_{q}(x)\ &=\ \lim_{N\rightarrow \infty}\sum_{x=0}^{p^{N}-1}f(x)\mu_{q}(x+p^{N}\mathbb{Z}_{p}) \label{1} \\
&=\ \lim_{N\rightarrow\infty}\frac{1}{[p^{N}]_{q}}\sum_{x=0}^{p^{N}-1}f(x)q^{x},\quad(\mathrm{see}\ [7,8]).\nonumber
\end{align}
From \eqref{1}, we have 
\begin{equation}
q\int_{\mathbb{Z}_{p}}f(x+1)d\mu_{q}(x)\ =\ \int_{\mathbb{Z}_{p}}f(x)d\mu_{q}(x)+(q-1)f(0)+\frac{q-1}{\log q}f^{\prime}(0), \label{2}
\end{equation}
where 
\begin{displaymath}
f^{\prime}(0)=\frac{df}{dx}\bigg|_{x=0},\quad(\mathrm{see}\ [1,2,7,8]). 
\end{displaymath}
Let us take $f(x)=e^{xt}$. Then, by \eqref{1}, we get 
\begin{equation}
\frac{(q-1)+\frac{q-1}{\log q}t}{qe^{t}-1}=\int_{\mathbb{Z}_{p}}e^{xt}d\mu_{q}(x). \label{3}
\end{equation}
The $q$-Bernoulli numbers are defined, in light of \eqref{3}, by
\begin{equation}
\frac{(q-1)+\frac{q-1}{\log q}t}{qe^{t}-1}=\sum_{n=0}^{\infty}B_{n,q}\frac{t^{n}}{n!}.\label{4}
\end{equation}
From \eqref{4}, we note that 
\begin{equation}
q(B_{q}+1)^{n}-B_{n,q}=\left\{\begin{array}{ccc}
	q-1, & \textrm{if $n=0$,} \\
	\frac{q-1}{\log q}, & \textrm{if $n=1$}, \\
	0, & \textrm{if $n > 1$}, 
\end{array}\right.\label{5}
\end{equation}
with the usual convention about replacing $B_{q}^{n}$ by $B_{n,q}$. \par 
For $|t|_{p}<p^{-\frac{1}{p-1}}$, the $(q,\lambda)$-Daehee polynomials are defined by 
\begin{equation}
\sum_{n=0}^{\infty}D_{n,q}(x|\lambda)\frac{t^{n}}{n!}=\frac{2(q-1)+\lambda\frac{q-1}{\log q}\log(1+t)}{q^{2}(1+t)^{\lambda}-1}(1+t)^{\lambda x},\quad(\mathrm{see}\ [3,12-17]). \label{6}	
\end{equation}
When $x=0$, $D_{n,q}(\lambda)=D_{n,q}(0|\lambda)$ are called $(q,\lambda)$-Daehee numbers.\par 
In particular, $\displaystyle D_{0,q}(0|1)=\frac{2}{[2]_{q}}\displaystyle$. \par 
The Catalan-Daehee numbers are defined by 
\begin{equation}
\frac{\frac{1}{2}\log(1-4t)}{\sqrt{1-4t}-1}=\sum_{n=0}^{\infty}d_{n}t^{n},\quad(\mathrm{see}\ [5,10]).\label{7}
\end{equation}
We note that 
\begin{equation}
\sqrt{1+t}=\sum_{m=0}^{\infty}(-1)^{m-1}\binom{2m}{m}\bigg(\frac{1}{4}\bigg)^{m}\bigg(\frac{1}{2m-1}\bigg)t^{m}.\label{8}	
\end{equation}
By replacing $t$ by $-4t$ in \eqref{8}, we get 
\begin{equation}
\sqrt{1-4t}=1-2\sum_{m=0}^{\infty}\binom{2m}{m}\frac{1}{m+1}t^{m+1}=1-2\sum_{m=0}^{\infty}C_{m}t^{m+1},\label{9}
\end{equation}
where $C_{m}$ is the Catalan number. \par 
From \eqref{7} and \eqref{9}, we have 
\begin{equation}
d_{n}=\left\{\begin{array}{ccc}
	1, & \textrm{if $n=0$,}\\
	\displaystyle\frac{4^{n}}{n+1}-\sum_{m=0}^{n-1}\frac{4^{n-m-1}}{n-m}C_{m}\displaystyle, & \textrm{if $n\ge 1$}.
\end{array}\right.\label{10}
\end{equation}
 When $q=1$, by \eqref{1}, we get 
 \begin{equation}
 \int_{\mathbb{Z}_{p}}(1-4t)^{\frac{x}{2}}d\mu_{1}(x)=\frac{\frac{1}{2}\log(1-4t)}{\sqrt{1-4t}-1}=\sum_{n=0}^{\infty}d_{n}t^{n}.\label{11}	
 \end{equation}

\section{$q$-analogues of Catalan-Daehee numbers and polynomials}
For $t\in\mathbb{C}_{p}$ with $|t|_{p}<p^{-\frac{1}{p-1}}$, we have 
\begin{equation}
\int_{\mathbb{Z}_{p}}(1-4t)^{\frac{x}{2}}d\mu_{q}(x)=\frac{q-1+\frac{q-1}{\log q}\frac{1}{2}\log (1-4t)}{q\sqrt{1-4t}-1}.\label{12}	
\end{equation}
In view of \eqref{11} and \eqref{12}, we define the $q$-analogue of Catalan-Daehee numbers by
\begin{equation}
\frac{q-1+\frac{q-1}{\log q}\frac{1}{2}\log (1-4t)}{q\sqrt{1-4t}-1}=\sum_{n=0}^{\infty}d_{n,q}t^{n}. \label{13}	
\end{equation}
Note that $\displaystyle\lim_{q\rightarrow 1}d_{n,q}=d_{n},\ (n\ge 0)\displaystyle$. \par 
From \eqref{6} and \eqref{13}, we have 
\begin{align}
\sum_{n=0}^{\infty}d_{n,q}t^{n}\ &=\ \frac{1}{2}\bigg(\frac{2(q-1)+\frac{q-1}{\log q}\log(1-4t)}{q^{2}(1-4t)-1}\bigg)\big(q\sqrt{1-4t}+1\big)\label{14} \\
&=\ \frac{1}{2}\bigg(\sum_{l=0}^{\infty}(-4)^{l}D_{l,q}(0|1)\frac{t^{l}}{l!}\bigg)\bigg(1+q-2q\sum_{m=0}^{\infty}C_{m}t^{m+1}\bigg) \nonumber\\
&=\ \frac{[2]_{q}}{2}\sum_{n=0}^{\infty}(-4)^{n}\frac{D_{n,q}(0|1)}{n!}t^{n}-q\sum_{n=1}^{\infty}\bigg(\sum_{m=0}^{n-1}\frac{(-4)^{n-m-1}}{(n-m-1)!}D_{n-m-1,q}(0|1)C_{m}\bigg)t^{n}\nonumber \\
&=\ 1+\sum_{n=1}^{\infty}\frac{[2]_{q}}{2}(-4)^{n}\frac{D_{n,q}(0|1)}{n!}t^{n}-q\sum_{n=1}^{\infty}\bigg(\sum_{m=0}^{n-1}\frac{(-4)^{n-m-1}}{(n-m-1)!}D_{n-m-1,q}(0|1)C_{m}\bigg)t^{n}\nonumber \\
&= 1+\sum_{n=1}^{\infty}\bigg(\frac{[2]_{q}}{2}\frac{(-4)^{n}}{n!}D_{n,q}(0|1)-q\sum_{m=0}^{n-1}\frac{(-4)^{n-m-1}}{(n-m-1)!}D_{n-m-1,q}(0|1)C_{m}\bigg)t^{n}. \nonumber
\end{align}
Therefore, by comparing the coefficients on both sides of \eqref{14}, we obtain the following theorem. 
\begin{theorem}
For $n\ge 0$, we have 
\begin{displaymath}
d_{n,q}=\left\{\begin{array}{ccc}
1, & \textrm{if $n=0$,}\\
\displaystyle\frac{[2]_{q}}{2}\frac{(-4)^{n}}{n!}D_{n,q}(0|1)-q\sum_{m=0}^{n-1}\frac{(-4)^{n-m-1}}{(n-m-1)!}D_{n-m-1,q}(0|1)C_{m}\displaystyle, & \textrm{if $n\ge 1$.}
\end{array}\right.
\end{displaymath}
\end{theorem}
From \eqref{13} and \eqref{14}, we have 
\begin{equation}
\sum_{n=0}^{\infty}\int_{\mathbb{Z}_{p}}x^{n}d\mu_{q}(x)\frac{t^{n}}{n!}=\int_{\mathbb{Z}_{p}}e^{xt}d\mu_{q}(x)=\frac{(q-1)+\frac{q-1}{\log q}t}{qe^{t}-1}=\sum_{n=0}^{\infty}B_{n,q}\frac{t^{n}}{n!}. \label{15}
\end{equation}
Thus, by \eqref{15}, we get 
\begin{equation}
\int_{\mathbb{Z}_{p}}x^{n}d\mu_{q}(x)=B_{n,q},\quad(n\ge 0).\label{16}	
\end{equation}
Now, we observe that 
\begin{align}
\sum_{n=0}^{\infty}d_{n,q}t^{n}\ &=\ \frac{q-1+\frac{q-1}{\log q}\frac{1}{2}\log (1-4t)}{q\sqrt{1-4t}-1}\ =\ \int_{\mathbb{Z}_{p}}(1-4t)^{\frac{x}{2}}d\mu_{q}(x)\label{17}\\
&=\ \sum_{m=0}^{\infty}\bigg(\frac{1}{2}\bigg)^{m}\frac{1}{m!}\big(\log(1-4t)\big)^{m}\int_{\mathbb{Z}_{p}}x^{m}d\mu_{q}(x)\nonumber\\
&=\ \sum_{m=0}^{\infty}\bigg(\frac{1}{2}\bigg)^{m}B_{m,q}\sum_{n=m}^{\infty}S_{1}(n,m)\frac{1}{n!}(-4t)^{n}\nonumber \\
&=\ \sum_{n=0}^{\infty}\bigg(\sum_{m=0}^{n}2^{2n-m}(-1)^{n}B_{m,q}S_{1}(n,m)\bigg)\frac{t^{n}}{n!}, \nonumber 
\end{align}
where $S_{1}(n,m),\ (n,m\ge 0)$ are the Stirling numbers of the first kind defined by 
\begin{displaymath}
(x)_{n}=\sum_{l=0}^{n}S_{1}(n,l)x^{l},\quad(n\ge 0),\quad(\mathrm{see}\ [1-17]).
\end{displaymath} 
Here $(x)_{0}=1$, $(x)_{n}=x(x-1)\cdots(x-n+1)$, $(n\ge 1)$. \par 
Therefore, by \eqref{17}, we obtain the following theorem. 
\begin{theorem}
For $n\ge 0$, we have 
\begin{displaymath}
(-1)^{n}d_{n,q}=\frac{1}{n!}\sum_{m=0}^{n}2^{2n-m}B_{m,q}S_{1}(n,m). 
\end{displaymath}
\end{theorem}
By binomial expansion, we get 
\begin{equation}
\int_{\mathbb{Z}_{p}}(1-4t)^{\frac{x}{2}}d\mu_{q}(x)=\sum_{n=0}^{\infty}(-4)^{n}\int_{\mathbb{Z}_{p}}\binom{\frac{x}{2}}{n}d\mu_{q}(x)t^{n}\label{18}.
\end{equation}
From \eqref{12}, \eqref{17} and \eqref{18}, we obtain the following corollary. 
\begin{corollary}
	For $n\ge 0$, we have 
	\begin{displaymath}
		\int_{\mathbb{Z}_{p}}\binom{\frac{x}{2}}{n}d\mu_{q}(x)\ =\ (-1)^{n}2^{-2n}d_{n,q}\ =\ \frac{1}{n!}\sum_{m=0}^{n}\bigg(\frac{1}{2}\bigg)^{m}B_{m,q}S_{1}(n,m).
	\end{displaymath}
\end{corollary}
The $q$-analogue of $\lambda$-Daehee polynomials are given by the following $p$-adic $q$-integral on $\mathbb{Z}_{p}$:
\begin{align}
\int_{\mathbb{Z}_{p}}(1+t)^{\lambda y+x}d\mu_{q}(y)\ &=\ \frac{(q-1)+\lambda\frac{q-1}{\log q}\log(1+t)}{q(1+t)^{\lambda}-1}(1+t)^{x} \label{19}\\
&=\ \sum_{n=0}^{\infty}D_{n,q,\lambda}(x)\frac{t^{n}}{n!}.\nonumber	
\end{align}
When $x=0$, $D_{n,q,\lambda}=D_{n,q,\lambda}(0),\ (n\ge 0)$, are called the $q$-analogue of $\lambda$-Daehee numbers. \par 
Here, we note that 
\begin{align}
\sum_{n=0}^{\infty}(-1)^{n}4^{n}D_{n,q,\frac{1}{2}}\frac{t^{n}}{n!}\ &=\ \frac{q-1+\frac{1}{2}\frac{q-1}{\log q}\log (1-4t)}{q(1-4t)^{\frac{1}{2}}-1} \label{20} \\
&=\ \sum_{n=0}^{\infty}d_{n,q}t^{n}. \nonumber 
\end{align}
Thus, by \eqref{20}, we get 
\begin{displaymath}
	d_{n,q}=(-1)^{n}\frac{4^{n}}{n!}D_{n,q,\frac{1}{2}},\quad (n\ge 0).
\end{displaymath}
Replacing $t$ by $\frac{1}{4}(1-e^{2t})$ in \eqref{13}, we have 
\begin{align}
\sum_{k=0}^{\infty}d_{k,q}\bigg(\frac{1}{4}\bigg)^{k}(1-e^{2t})^{k}\ &=\ \frac{q-1+\frac{q-1}{\log q}t}{qe^{t}-1}\ =\ \int_{\mathbb{Z}_{p}}e^{xt}d\mu_{q}(x)\label{21} \\
&=\ \sum_{n=0}^{\infty}B_{n,q}
\frac{t^{n}}{n!}. \nonumber	
\end{align}
On the other hand, 
\begin{align}
\sum_{k=0}^{\infty}d_{k,q}\bigg(\frac{1}{4}\bigg)^{k}(1-e^{2t})^{k}\ &=\ \sum_{k=0}^{\infty}k!d_{k,q}\bigg(-\frac{1}{4}\bigg)^{k}\frac{1}{k!}\big(e^{2t}-1\big)^{k} \label{22}\\
&=\ \sum_{k=0}^{\infty}k!d_{k,q}\bigg(-\frac{1}{4}\bigg)^{k}\sum_{n=k}^{\infty}S_{2}(n,k)2^{n}\frac{t^{n}}{n!}\nonumber\\
&=\ \sum_{n=0}^{\infty}\bigg(\sum_{k=0}^{n}(-1)^{k}k!d_{k,q}2^{n-2k}S_{2}(n,k)\bigg)\frac{t^{n}}{n!},\nonumber 
\end{align}
where $S_{2}(n,k),\ (n,k\ge 0)$, are the Stirling numbers of the second kind defined by 
\begin{displaymath}
	x^{n}=\sum_{l=0}^{n}S_{2}(n,l)(x)_{l},\quad(n\ge 0). 
\end{displaymath}
Therefore, by \eqref{21} and \eqref{22}, we obtain the following theorem. 
\begin{theorem}
For $n\ge 0$, we have 
\begin{displaymath}
B_{n,q}=\sum_{k=0}^{n}(-1)^{k}2^{n-2k}k!S_{2}(n,k)d_{k,q}. 
\end{displaymath}
\end{theorem}
Now, we observe that 
\begin{displaymath}
	\int_{\mathbb{Z}_{p}}(1-4t)^{\frac{x+y}{2}}d\mu_{q}(y)=\frac{(q-1)+\frac{q-1}{\log q}\frac{1}{2}\log(1-4t)}{q\sqrt{1-4t}-1}(1-4t)^{\frac{x}{2}}.
\end{displaymath}
We define the Catalan-Daehee polynomials by 
\begin{equation}
\frac{q-1+\frac{q-1}{\log q}\frac{1}{2}\log(1-4t)}{q\sqrt{1-4t}-1}(1-4t)^{\frac{x}{2}}=\sum_{n=0}^{\infty}d_{n,q}(x)t^{n}. \label{23}
\end{equation}
Note that 
\begin{align}
(1-4t)^{\frac{x}{2}}\ &=\ \sum_{l=0}^{\infty}\bigg(\frac{x}{2}\bigg)^{l}\frac{1}{l!}\big(\log(1-4t)\big)^{l}\ =\ \sum_{l=0}^{\infty}\bigg(\frac{x}{2}\bigg)^{l}\sum_{m=l}^{\infty}S_{1}(m,l)(-4)^{m}\frac{t^{m}}{m!} \label{24}\\
&=\ \sum_{m=0}^{\infty}\bigg(\sum_{l=0}^{m}S_{1}(m,l)\frac{(-4)^{m}}{m!}\bigg(\frac{x}{2}\bigg)^{l}\bigg)t^{m}.\nonumber	
\end{align}
Thus, by \eqref{13}, \eqref{23} and \eqref{24}, we get 
\begin{align}
\sum_{n=0}^{\infty}d_{n,q}(x)t^{n}\ &= \frac{q-1+\frac{q-1}{\log q}\frac{1}{2}\log(1-4t)}{q\sqrt{1-4t}-1}(1-4t)^{\frac{x}{2}}\label{25} \\
&=\ \sum_{k=0}^{\infty}d_{k,q}t^{k}\sum_{m=0}^{\infty}\bigg(\sum_{l=0}^{m}S_{1}(m,l)\frac{(-4)^{m}}{m!}\bigg(\frac{x}{2}\bigg)^{l}\bigg)t^{m}\nonumber 	\\
&=\ \sum_{n=0}^{\infty}\bigg(\sum_{m=0}^{n}\sum_{l=0}^{m}S_{1}(m,l)\frac{(-4)^{m}}{m!}d_{n-m,q}\bigg(\frac{x}{2}\bigg)^{l}\bigg)t^{n}. \nonumber
\end{align}
By comparing the coefficients on both sides \eqref{25}, we obtain the following theorem. 
\begin{theorem}
For $n\ge 0$, we have 
\begin{equation*}
d_{n,q}(x)\ =\  \sum_{l=0}^{n}\bigg(\sum_{m=l}^{n}(-1)^{m}\frac{2^{2m-l}}{m!}S_{1}(m,l)d_{n-m,q}\bigg)x^{l}. 
\end{equation*}
\end{theorem}

\section{Conclusion}

Quite a few special numbers and polynomials have been studied by employing various different tools. Previously, the Catalan-Daehee numbers and polynomials were introduced by means of $p$-adic Volkenborn integrals and some interesting results for them were obtained by using generating functions, differential equations, umbral calculus and $p$-adic Volkenborn integrals. In this paper, we introduced $q$-analogues of the Catalan-Daehee numbers and polynomials and obtained several explicit expressions and identities related to them. In more detail, we expressed the Catalan-Daehee numbers in terms of the $(q,\lambda)$-Daehee numbers, and of the $q$-Bernoulli polynomials and Stirling numbers of the first kind. We obtained an identity involving $q$-Bernoulli number, $q$-analogues of Catalan-Daehee numbers and Stirling numbers of the second kind. In addition, we got an explicit expression for the $q$-analogues of  Catalan-Daehee polynomials which involve the $q$-analogues of Catalan-Daehee numbers and Stirling numbers of the first kind. \par
It has been our constant interest to find $q$-analogues of some interesting special numbers and polynomials and to study their arithmetic and combinatorial properties and their applications. We would like to continue to study this line of research in the future.

\end{document}